\documentclass{article}
\usepackage{enumitem}
\title{A Mathematician Reads the Kalam Cosmological Argument}
\author{Timothy Y. Chow}
\date{July 2022}

\begin{document}
\maketitle

\section{Introduction}

When I was an undergraduate, I had an interest in Christian
apologetics, and one of the books I studied was by
William Lane Craig~\cite{craig-apologetics}.  Reading
through the early chapters, I felt that it was one of
the better books of the genre---until I came to the
section on something Craig called the \emph{kalam
cosmological argument} for the existence of God.
When I read the argument, I cringed. Not only did
the argument strike me as fallacious; it seemed
\emph{obviously} fallacious, or at least obviously
fallacious to anyone with a certain amount of background in
mathematics and physics.  I finished
reading the book anyway, but was left feeling that the
section on the kalam cosmological argument was a huge flaw
in an otherwise decent book.

Recently, the kalam cosmological argument came to
my attention again, and out of curiosity,
I searched the literature for critiques.
I found many, but to my surprise,
nobody seemed to have articulated the objection
that had seemed so obvious to my younger self.
I decided to write out my objection in some detail,
and I shared the draft with about a dozen people
with some knowledge of the subject.
What I found striking was that
the best predictor of whether someone would agree
with me was not their religious beliefs,
but their academic background.
Physicists and mathematicians generally
agreed with me and did not think
I was saying anything particularly controversial,
whereas philosophers would balk and raise what I considered
to be peculiar counter-objections.
There seems to be some kind of culture gap or language gap;
the hard-earned wisdom acquired by the mathematical
and scientific communities through decades if not centuries
of experience is failing to be transmitted to the
philosophical community.
I find this culture gap troubling,
so I am publishing this essay with the
hope that I can help bridge it.

But I am getting ahead of myself. Just what is the kalam
cosmological argument, anyway?

\section{The Kalam Cosmological Argument}

In apologetics, the term
\emph{cosmological argument},
or argument from first cause, refers to an argument
with something like the following structure.
\begin{enumerate}
\item Whatever begins to exist has a cause of its beginning.
\item The universe began to exist.
\item Therefore, the universe has a cause of its beginning.
\end{enumerate}
Since 1979, William Lane Craig (as well as other Christian apologists, but for brevity I will use Craig’s name synecdochically throughout this article) has vigorously championed a version of the cosmological argument called the \emph{kalam cosmological argument}.
The literature on the kalam cosmological argument is remarkably large and continues to grow, with several entire books, some quite recent, devoted to it~\cite{copan-craig-1,copan-craig-2,erasmus,nowacki,pearce-lowder-east}.

The kalam cosmological argument itself has many
sub-variants,
but in this paper, I will focus on a specific version,
which argues for Claim~2 above
(the universe began to exist) as follows.
\begin{enumerate}[label=\Alph*.]
\item An actual infinite cannot exist
(in the physical world).
\item An infinite temporal regress is an actual infinite.
\item Therefore, an infinite temporal regress of events cannot exist.
\end{enumerate}
Even more specifically, I want to criticize Craig's
argument for sub-claim~A, that
\emph{an actual infinite cannot exist}
in the physical world.
This sub-claim is perhaps the feature that most clearly
distinguishes the kalam cosmological argument from
other cosmological arguments.

We will examine some of Craig's arguments in detail
later, but for now, it is important to notice that
Craig does not just argue that an actual infinite
\emph{does not} exist; he argues that an actual
infinite \emph{cannot} exist.
In other words, its existence is \emph{impossible}.
Furthermore, the sense of impossibility that Craig
intends here is \emph{metaphysical impossibility}.
That is, according to Craig,
the problem with the existence of an actual
infinite is not that it is \emph{logically} absurd
or that it contradicts experimental observations;
the problem is that it violates some metaphysical
principle that Craig considers to be inviolable.

The astute reader may have already guessed where I am headed.
Scientists know from bitter experience that insisting
that the physical world must conform to some
metaphysical principle that our inner voice tells us
must be true, even when no logical argument or
empirical evidence requires it, has been shown
time and time again to be a rash and highly unreliable
approach to physics.  This, in a nutshell, is why
I find Craig's argument that an actual infinite
\emph{cannot} exist to be completely unconvincing.

The plan of this essay is as follows.
First, as a preliminary step, I will clarify different
senses of impossibility---logical, physical, and
metaphysical.  I will then recall some famous
historical examples of failed attempts to argue on
the basis of metaphysical impossibility;
the scientific community has learned the hard way
to shun this type of metaphysical argument.
After laying this groundwork, I will examine
some of Craig's arguments specifically.
It will, I hope, become apparent that his supposedly
inviolable metaphysical principles are even more
flimsy and unconvincing than the principles that
great thinkers of the past, including Einstein,
mistakenly thought were mandatory.

\section{Logical and Physical Impossibility}

As mentioned already, Craig is insistent
that the existence of an actual infinite is not just
\emph{false} but is \emph{impossible}.
But what does it mean to say that something is
``impossible''?

One way that something can be impossible is for it to be
\emph{logically} impossible or \emph{mathematically}
impossible\footnote{There is a philosophy of mathematics
known as \emph{logicism} that equates mathematics with logic.
I do not necessarily endorse logicism, but I believe that
the distinction between logical possibility and mathematical
possibility does not matter for the present discussion.}.
For example, one could claim that the very concept of an
actual infinite---that is, a completed infinite totality of
objects---is incoherent or logically self-contradictory.
In particular, one could claim that Zermelo–Fraenkel set
theory (ZF for short, or ZFC if one includes the axiom of
choice), a widely accepted axiomatization of most of modern
mathematics that among other things postulates the existence
of an infinite set, is inconsistent.

Craig has made it clear that he is not claiming that an
actual infinite is logically impossible, and in particular,
he accepts that ZF is consistent.  In this regard, Craig is
more generous than some other critics of the actual infinite.
Until the work of Georg Cantor in the 19th century,
most thinkers accepted only
something they called \emph{potential infinity},
and rejected actual infinities as being conceptually
suspect and maybe even logically incoherent\footnote{By
contrast, Cantor himself believed not only
that actual infinities existed abstractly, but
``in the natura creata an Actual Infinite must be
confirmed, as for example with respect to,
in my strong conviction, the actual infinite number of
created individual beings, not only in the universe
but also already on our earth''~\cite{cantor}.}.
Even today, there are skeptics who do not
take for granted that ZF is consistent.
Such doubts are interesting in their own right,
but we will set them aside,
because Craig does not raise them.
He accepts that actual infinities are conceptually
coherent, and is concerned only with their
physical existence.  After all,
even if something is logically possible, it might be
\emph{physically} impossible.

It is harder to define
physical impossibility precisely than it is to define
mathematical impossibility precisely, but here is a
definition that should suffice for our present purposes.
Let us use the term \emph{physical theory} for a mathematical
model of the physical world that makes empirically testable
predictions about the physical world.
Let us then say that a physically observable event or
circumstance---call it $X$---is physically impossible if our
best physical theories predict that $X$ will never happen.
For example, let $X$ be the transmission of information
faster than the speed of light.  $X$ is logically possible
and there is clear sense in which we could observe $X$
happening---for example, in 2011, there was a high-profile
experiment that at first (before flaws were discovered in
the experimental setup) seemed to show that neutrinos
traveled faster than light from CERN to the Gran Sasso
National Laboratory in Italy.
However, our best mathematical models of the physical world
predict that information cannot travel faster than the speed
of light.  Similarly, our best physical theories of
statistical mechanics and thermodynamics predict that any
machine we build will never exhibit perpetual motion---or
at least that the probability of its doing so is absurdly
small.

Notice carefully a couple of features of this definition of
physical impossibility.  First, a clear distinction is made
between our mathematical models of the physical world and
the physical world itself.  This distinction is a hallmark
of all modern science.  In particular, physical impossibility
as we have defined it is not directly predicated of the
physical world itself.  Arguably, the real world just is the
way it is, and if so, it is not totally clear what it would
mean for something that is \emph{not true} in the real world
to be \emph{possible}.  On the other hand, in a mathematical
model of the physical world, we can hypothesize a variety of
circumstances and give a clear meaning to possibility and
impossibility.

Second, a clear distinction is made between the mathematical
model and its physical predictions.
A mathematical model might have all sorts of exotic
mathematical features---including the presence of infinite
sets---but not all such features necessarily yield observable
predictions, and they may or may not correspond directly to
features of the physical world.
It is only the model's predictions about the physical world
that necessarily have physical meaning.  In particular, I
would defend Craig against one of his critics,
Quentin Smith~\cite{smith}, who says that ``Big bang
cosmology implies that there is an actually infinite
manifold, topology, and metrication.''
The usual \emph{mathematical model}
of big bang cosmology does indeed
contain such actually infinite things, but even if we accept
the model, we are not automatically committing ourselves to
asserting that the \emph{physical world} contains actual
infinities.  Craig's arguments against an actual infinite
cannot be refuted so cheaply.

\section{Metaphysical Impossibility?}

Craig certainly believes that an actual infinite is
physically impossible.  However, as already mentioned,
Craig argues for something stronger---namely that an actual
infinite is \emph{metaphysically} impossible,
and it is this aspect of Craig's arguments that
I take issue with. But what does
it mean for something to be metaphysically impossible?

The exact meaning of metaphysical impossibility is a
notoriously controversial philosophical issue,
but roughly speaking,
it means a kind of impossibility that is not sensitively
dependent on accidental or contingent facts about our
current mathematical models or the physical universe we
happen to live in, but holds for more general reasons.
At minimum, if something is metaphysically impossible,
then not only must it be impossible according to the
physical theories that we currently deem best;
it must also be impossible according to any possible
physical theory that we might someday deem best,
and it must also be physically impossible in conceptually
coherent hypothetical universes other than the actual
universe.

Suppose we grant that metaphysical impossibility
is a meaningful concept.
Then the following 
strategy for arguing that $X$ is metaphysically
impossible---and a fortiori physically
impossible---is tempting and, as we shall see later,
is exactly how Craig argues.
\begin{enumerate}
\item Decide a priori on some metaphysical principle $P$
that one feels must be obeyed
by any satisfactory physical theory.
\item Show that $X$ is incompatible with $P$, using abstract
argumentation only and not the results of any physical
observation or experiment.
\item Conclude that $X$ is absurd or impossible.
\end{enumerate}
Unfortunately, as many historical examples illustrate,
the above strategy has an embarrassingly bad track record.

For our first example, let us recall Giovanni Girolamo
Saccheri, a Jesuit priest who, in the 18th century,
published what is now regarded as a mathematical
investigation into non-Euclidean geometry~\cite{saccheri}.
Saccheri considered what would happen if we were to reject
the fifth postulate of Euclidean geometry and instead were
to suppose that the sum of the angles of a triangle were
less than 180 degrees.
Saccheri failed to deduce any mathematical contradiction,
yet he famously declared in Proposition XXXIII,
``\emph{Hypothesis anguli acuti est absoluta falsa;
quia repugnans naturae lineae rectae}.''
(The acute angle hypothesis is absolutely false,
because it is repugnant to the nature of the straight line.)

It is important to remember that, at the time, Euclidean geometry was regarded not just as a purely imaginary logical construction, but as a mathematical model of physical space.
Therefore, Saccheri was in effect declaring that the phenomena of non-Euclidean geometry not only \emph{were} not,
but \emph{could not be}, part of a valid physical theory.
There is some debate about whether Saccheri, in his heart of hearts, truly believed that non-Euclidean geometry was absurd, or whether he was saying such things to protect himself from anticipated criticism,
but either way, we have here a clear example of a declaration of impossibility based purely on a metaphysical prejudice about what ``the nature of a straight line'' must be, and not on any logical inconsistency or experimental tests.
Today, of course, not only do we know that
non-Euclidean geometry is perfectly consistent mathematically, but we have good reasons to believe that the
actual geometry of spacetime is non-Euclidean.

Our second example concerns an old and (at one time) widely believed metaphysical principle that if something, such as light, exhibits wave behavior,
then it \emph{must} be supported by some kind of medium.
Indeed, it seemed part of the very concept of a wave that it was a disturbance in an underlying medium.
In the case of light, it was hypothesized, on the basis of this metaphysical principle and not on any direct experimental evidence,
that there existed a medium called the
\emph{luminiferous aether}.
But after the famous Michelson--Morley experiment and the development of Einstein’s theory of relativity,
the aether hypothesis was eventually discarded, and with it, the claim that it was impossible for a wave to propagate through a vacuum.

For our third example, let us recall that
Albert Einstein~\cite{einstein}
famously rejected quantum mechanics, saying,
``\emph{Die Quantenmechanik ist sehr achtunggebietend.
Aber eine innere Stimme sagt mir, da\ss\
das noch nicht der wahre Jakob ist.
Die Theorie liefert viel, aber dem Geheimnis des Alten
bringt sie uns kaum n\"{a}her.
Jedenfalls bin ich \"{u}berzeugt,
da\ss\ der nicht w\"{u}rfelt}.”
(Quantum mechanics is very impressive.  But an inner voice tells me that it is not yet the real McCoy\footnote{I have borrowed the colorful translation
``the real McCoy'' from Ryckman~\cite{ryckman},
who also suggests that \emph{der wahre Jakob},
in addition to being a biblical reference,
may be an allusion to a satirical German newspaper
of that name.}.
The theory yields much, but it hardly brings us closer to the Ancient One's secrets.  I, at any rate, am convinced that he does not throw dice.)
Along with Podolsky and
Rosen~\cite{einstein-podolsky-rosen},
Einstein proposed a thought experiment that we would nowadays describe as exhibiting \emph{quantum entanglement}.
For example, one can prepare an
entangled electron-positron pair\footnote{The observables
in the original EPR paper were position and momentum,
but the phenomenon is easier to describe using discrete
observables such as spin.}
so that (1)~each particle’s spin when measured along (say) the $x$-axis is equally likely to be
``spin up'' or ``spin down'' and
(2)~when measured along the same axis,
the positron's spin is always the opposite of the electron's
spin, even if the two particles are separated so far apart that they cannot communicate their ``decisions''
to each other without violating the speed-of-light barrier.  This state of affairs is unsurprising if there is
a ``hidden variable'' that determines the state of the system
before the particles are separated and measured, but seems
absurd if the spins are indeterminate.
Einstein, being committed to a metaphysical principle of
realism that required quantum-mechanical observables such as
spin to have a definite value at all times, felt that this
thought experiment showed that quantum mechanics had to be
incomplete.  However, later work by Bell, Aspect, Kochen,
Specker, and others has shown that quantum entanglement
is real, and that the kind of
local hidden variable theory that Einstein hoped for is
inconsistent with known experimental facts\footnote{For
a thorough discussion of this point, see Wiseman
and Cavalcanti~\cite{wiseman}, and for experimental
evidence, see Hensen et~al.~\cite{hensen}}.
It is Einstein's metaphysical principle that has been shown
to be dubious.

It would be easy to give more examples from modern physics
that might seem metaphysically
impossible---light behaving like a wave as well as a
particle, identical twins aging at different rates simply
because one of the twins has made a very speedy round trip
to a distant location (the \emph{twin paradox} of relativity
theory), the universe possibly having ten dimensions,
and so on.  Physicists have learned the hard way
to obey what we might dub a ``prudence principle'':

\begin{quote}
\emph{If a physical theory is consistent with known
experimental and observational facts, then do not reject
it out of hand purely because it contradicts some
intuitively plausible metaphysical principle.}
\end{quote}
After all, if Einstein’s intuition was fallible,
who dares claim to be immune?

If you propose a crazy new theory, scientists will naturally
react skeptically, but if the only thing your theory
contradicts is some metaphysical principle and not any
mathematical calculation or experimental fact, then they
will probably not dismiss the theory as impossible
\emph{a priori}.
Instead, their instinct will be to ask, can you design an
experiment to test your theory, and mathematically calculate
the predicted outcome?  Whenever possible, scientists
gravitate towards mathematics and/or experiment and not
metaphysics.

The prudence principle embodies a certain kind of skepticism
of metaphysical dogmatism, and may remind the reader of
logical positivism or of Karl Popper’s famous concept of
\emph{falsifiability}, but there is a crucial difference.
For example, suppose you have a favorite metaphysical
principle~$P$.  Popper would reject $P$
as unscientific unless
you could articulate some conceivable experimental outcome
that would refute it.  On the other hand, the prudence
principle makes the weaker claim that you should not insist
on $P$ to the point where you declare the alternative to be
absurd.  While falsifiability is widely accepted by
scientists, it is not entirely uncontroversial; for example,
physicist Sean Carroll has criticized it.  However, even
Carroll accepts (a version of) the prudence principle. 
Regarding certain theories that are arguably unfalsifiable,
Carroll~\cite{carroll} has written,
``Refusing to contemplate their
possible existence on the grounds of some a priori principle,
even though they might play a crucial role in how the world
works, is as non-scientific as it gets.''

Far from a blanket rejection of metaphysics, the prudence
principle is a call for metaphysical tolerance, and not a
call to positivism or scientism.  It should of course not
be thought of as an ironclad rule, but it represents
hard-won wisdom, and I maintain that a heavy burden of
proof rests on anyone who wants to violate~it.
A metaphysical principle had better be \emph{extremely}
compelling if we are going to use it to reject
physical theories a priori.
As we shall soon see, none of Craig's metaphysical
principles comes remotely close to qualifying.

\section{Hilbert's Hotel}

Let us now turn to one of Craig's favorite arguments against
the existence of an actual infinite, namely Hilbert's hotel,
a thought experiment due to the German mathematician David
Hilbert.  Most readers of the \emph{Intelligencer} will be
familiar with Hilbert's Hotel, but for the sake of
completeness, let us review the setup.
We imagine a hotel with a countably infinite number of
rooms, labeled 1, 2, 3, etc.  Suppose that every room is
occupied.  Now a countably infinite number of new guests
show up.  The proprietor accommodates them all by shifting
the person in room $n$ to room $2n$ for all $n$, thereby
vacating infinitely many rooms for the new guests.
Craig and Sinclair~\cite{craig-sinclair}
continue the narrative as follows.

\begin{quote}
But Hilbert's Hotel is even stranger than the German
mathematician made it out to be.  For suppose some of the
guests start to check out. Suppose the guest in room \#1
departs.  Is there not now one fewer person in the hotel?
Not according to infinite set theory!  Suppose the guests
in rooms \#1, 3, 5, $\ldots\,$ check out. In this case an
infinite number of people has left the hotel, but by Hume's
Principle, there are no fewer people in the hotel. In fact,
we could have every other guest check out of the hotel and
repeat this process infinitely many times, and yet there
would never be any fewer people in the hotel. Now suppose
the proprietor does not like having a half-empty hotel (it
looks bad for business). No matter!  By shifting guests in
even-numbered rooms into rooms with numbers half their
respective room numbers, he transforms his half-vacant
hotel into one that is completely full. In fact, if the
manager wanted double occupancy in each room, he would
have no need of additional guests at all. Just carry out
the dividing procedure when there is one guest in every
room of the hotel, then do it again, and finally have one
of the guests in each odd-numbered room walk next door to
the higher even-numbered room, and one winds up with two
people in every room!

One might think that by means of these maneuvers the
proprietor could always keep this strange hotel fully
occupied. But one would be wrong. For suppose that the
persons in rooms \#4, 5, 6, $\ldots\,$ checked out.
At a single stroke the hotel would be virtually emptied,
the guest register reduced to three names, and the infinite
converted to finitude.  And yet it would remain true that
as many guests checked out this time as when the guests in
rooms \#1, 3, 5, $\ldots\,$ checked out! Can anyone believe
that such a hotel could exist in reality?

Hilbert's Hotel is absurd. But if an actual infinite were
metaphysically possible, then such a hotel would be
metaphysically possible. It follows that the real existence
of an actual infinite is not metaphysically possible.
\end{quote}

Again, most readers of the \emph{Intelligencer} probably
do not need me to point out the flaws in the above
reasoning, but let us go ahead and spell some of them out.
First of all,
as other critics have noted, the claim that ``if an actual
infinite were metaphysically possible, then such a hotel
would be metaphysically possible'' is obscure.  By analogy
with the concept of NP-completeness in computational
complexity theory, we might rephrase Craig and Sinclair's
claim as a claim that Hilbert's hotel is
\emph{AI-complete}
(where ``AI'' stands for ``Actual Infinite'');
i.e., that the impossibility of an actual infinite reduces
to the impossibility of Hilbert's hotel.
But why is Hilbert's hotel AI-complete?
Craig and Sinclair do not explain this point.
In fact, it looks as though they commit an error
that students of complexity theory often make,
which is to get the direction of the reduction backward.
If all actual infinities are impossible, then
in particular, Hilbert's hotel is impossible;
that much is obvious.
But the converse is not at all obvious.

But let us not dwell on this objection,
since a more important
question is, what justifies the conclusion that Hilbert's
hotel is absurd?

There does seem to be something absurd about having
one person leave and yet not having one fewer person.
But this absurdity comes not from infinite set theory,
as Craig and Sinclair claim, but from a linguistic quibble.
If by \emph{fewer} we mean \emph{strictly lower
cardinality} then indeed, removing a single element from
an infinite set does not lower its cardinality,
but there is nothing absurd about that.
On the other hand,
if $A$ and $B$ are sets and $A\subseteq B$,
then we could choose to define
``$A$ has one fewer member than~$B$''
by $|B\backslash A| = 1$, where $\backslash$
denotes set subtraction.  With that definition,
the set of guests after the departure does indeed
have one fewer member than the set of guests before
the departure.  Any appearance of absurdity evaporates
as soon as one clearly defines one's terms.

The other observations that Craig and Sinclair make about
Hilbert's hotel may seem strange to some, but where is the
absurdity?  To those who work with infinite sets on a daily
basis, everything in the account seems normal and
nothing in sight seems absurd.

As mentioned earlier, Craig and Sinclair try to
state more explicitly what they think is absurd by
stating metaphysical principles that they regard as
inviolate.
For example, 
one metaphysical principle $P$ that Craig and
Sinclair propose is that physical quantities
(i.e., numbers of physically existing entities)
must obey the law that equal quantities can always be
subtracted from equal quantities and the resulting
quantities must always be equal.  There are different
ways to formulate $P$ precisely; one possibility is this:
\begin{quote}
$P$: If $A_1\subseteq B_1$ and $A_2 \subseteq B_2$
and $ |A_1| = |A_2|$ and $|B_1| = |B_2|$ then
$|B_1| - |A_1|$ and $|B_2| - |A_2|$
always exist and are equal to each other
and to $|B_1\backslash A_1|$ and $|B_2\backslash A_2|$.
\end{quote}
It is indeed true that $P$ does not necessarily hold for
infinite sets.  Craig and Sinclair not only maintain
that $P$ is plausibly true for physical quantities;
they make the far stronger claim that $\neg P$,
the negation of~$P$, is \emph{absurd}.
But why?  Why should we believe in~$P$ so strongly for
physical quantities that we should rule out a priori any
physical theory that violates it, even when no logical
principle or experimental fact is contradicted?

My objection here is not entirely new.
Other critics have pointed out that $\neg P$
just comes with the territory of infinite sets
and is not the absurdity that Craig claims it is.
Craig's reaction has been to say that his critics
do nothing to address the apparent absurdity of $\neg P$
for physical quantities.  But Craig misplaces the burden
of proof.  The burden of proof is rather on Craig to
demonstrate that $\neg P$ is so repugnant to the nature
of the universe that it justifies throwing caution, in the
form of the prudence principle, to the winds.

Another metaphysical principle suggested by Craig and
Sinclair is that ``it is ontologically absurd that a hotel
exist which is completely full and yet can accommodate
untold infinities of new guests just by moving people
around.''  Once
again we are led to ask, what exactly is so absurd
about it?  It seems far easier to swallow than
quantum entanglement or the twin paradox.  The only
argument that Craig and Sinclair present is
an argument from incredulity.

To make matters worse, it is not even clear that Hilbert's
hotel is \emph{physically} impossible, let alone
metaphysically impossible.  Of course, we human beings are
finite creatures and could not possibly build Hilbert's
hotel, nor could we engineer infinite evacuations and
reassignments.  However, the question is not whether we
humans have the capacity to carry out infinite operations,
but whether Hilbert's hotel could exist.  To put it another
way, if we were to send out a space probe and it were to
stumble upon what looked like the beginning of Hilbert's
hotel floating in interstellar space, is there anything in
our current best physical theories that would demand that
the hotel be finitely long?

The answer seems to be no.  Now, if it were the case
that our current best physical
theories predicted that the universe were spatially finite
(e.g., that its spatial topology were like that of a
three-dimensional torus), then they would also predict
that any sufficiently long (non-self-intersecting) hotel
would fill up the universe.  In particular, Hilbert's
hotel would be physically impossible.  However, the current
scientific consensus is that the spatial finitude of the
universe is an open question\footnote{See for example
Levin et~al.~\cite{levin}
More recent calculations in a similar vein have
yielded the same qualitative conclusion.}, so the physical
impossibility of Hilbert's hotel cannot be deduced in
this manner.  Let us note also that if a spatially
infinite universe were really as absurd as Craig and
Sinclair claim, it is curious that astrophysicists
seem not to have noticed this absurdity and continue to
regard it as a viable theory.

\section{Benardete's Book}

Here is a thought experiment by Jos\'e
Benardete~\cite{benardete} that Craig likes
to cite as an argument against an actual infinite.

\begin{quote}
Here is a book lying on the table.  Open it.  Look at the
first page.  Measure its thickness.  It is very thick indeed
for a sheet of paper---1/2 inch thick.  Now turn to the
second
page of the book.  How thick is this second sheet of paper?
1/4 inch thick.  And the third page of the book, how thick
is this third sheet of paper?  1/8 inch thick, \&c.\
\emph{ad infinitum}.
We are to posit not only that each page of the
book is followed by an immediate successor the thickness of
which is one-half that of the immediately preceding page but
also (and this is not unimportant) that each page is
separated from page 1 by a finite number of pages.  These
two conditions are logically compatible: there is no
certifiable contradiction in their joint assertion.
But they mutually entail that there is no last page in the
book.  Close the book.  Turn it over so that the front
cover of the book is now lying face down upon the table.
Now---slowly---lift the back cover of the book with the
aim of exposing to view the stack of pages lying beneath
it.  \emph{There is nothing to see}.
For there is no last page
in the book to meet our gaze. [Emphasis in original.]
\end{quote}

Benardete's book has a feature that Hilbert's hotel lacks.
Namely, there is an obvious reason why Benardete's book is
\emph{physically} impossible: According to our best physical
theories, in particular the so-called Standard Model of
particle physics, a sheet of paper cannot be arbitrarily
thin.  However, remember that Craig wishes to argue for
more than the physical impossibility of Benardete's book;
he wants to argue that its impossibility does not depend
on accidental features of our current best physical
theories, but is a general fact that must hold of all
possible physical theories that we might someday deem best.

But it is totally unclear why Benardete's book could not
exist in some hypothetical physical universe that is similar
to ours in many respects but in which pages of a book can
have arbitrarily small finite thickness.
At first glance, it does seem absurd that there would
be ``nothing to see'' if we were
to open the back of the book.
But let us think more carefully about the situation.
In the mathematical model
of Benardete's book, there is indeed no last page,
but in order to make a physical
prediction, we must flesh out some further details.
Vision, in our physical universe, is
based on electromagnetic radiation in the visible part of
the spectrum.  Light scatters off a page and some of it
enters our eyes.  In order to predict what we would see
when we lift the back cover of Benardete's book, we need
some model of vision in this hypothetical universe.  Many
possibilities suggest themselves.  Perhaps in this universe,
there are waves that are somewhat similar in character to
light waves, but there is some finite thickness $d$ such
that the light waves pass through pages with thickness less
than $d$ without interacting, so that what we would see
would be the last page whose thickness is at least~$d$.
There is nothing absurd in this scenario;
there would be infinitely many invisible pages,
with no last page, but it would not be the case
that there would be ``nothing to see.''

To get a contradiction, we have to postulate some
metaphysical principle---perhaps a claim that the only
thing that can cause a light wave to scatter is a specific
page, and it cannot scatter off a specific page if there
is another page in front of it.  But what forces us to
adopt this metaphysical principle?
Is the argument just that the alternative is repugnant to
the nature of a page?  Benardete's book is
no more bizarre than non-Euclidean geometry or a wave in
a vacuum, and we are offered no credible reason to
reject it a priori.
Ironically, in his essay, Benardete himself argues strongly
\emph{against} Craig's belief
that we know a priori, apart from all empirical
evidence, that something like Hilbert's hotel cannot exist.

\section{Nowacki's Substance-Based Metaphysics}

There are various other thought experiments that Craig has
cited, but they all have a similar flavor, and in every case
the weakness is the same---a contradiction is claimed, but
the only thing being contradicted is some metaphysical
postulate that Craig finds intuitively obvious; e.g., that
the number of orbits that a planet makes must be either
even or odd (and therefore cannot be infinite).
Insisting on such metaphysical postulates
simply because the alternative seems weird
exhibits dogmatism of
a type that the entire scientific community knows from
long experience to avoid.

I said earlier that when I have presented my argument
to philosophers, they often balk.
I have found that appealing to principles
(such as what I have been calling the prudence principle)
that are backed primarily by the
collective experience of subject-matter experts,
rather than by deductive reasoning,
is typically not recognized by philosophers as
a convincing mode of argumentation.
More acceptable to them is to propose a metaphysical
principle
that seems plausible to philosophers,
no matter how unfounded it may seem to
expert practitioners of the subject in question
(in this case, physicists).
This seems to be a cultural gap between the
two academic communities that is difficult to bridge.

To be fair, however, not every philosopher who
argues for the impossibility of the actual infinite
pulls arbitrary metaphysical principles out of thin air
and expects others to accept them as obvious.
For example, Nowacki~\cite{nowacki}
is sympathetic to the kalam cosmological argument,
but he correctly recognizes that one must carefully
articulate and defend the metaphysical assumptions that are
being contradicted.
I applaud Nowacki for recognizing this point.
Unfortunately, I have serious doubts about his suggested
metaphysical framework.  Nowacki proposes what he calls a
\emph{substance-based metaphysics}.
It is not entirely clear exactly what Nowacki means by a
\emph{substance}, but he does seem to agree with what he
calls the commonsense position of Aristotle, ``who defined
individual substance as what exists without either being
predicated of or existing in anything else.''
Moreover, Nowacki intends substances to be medium-sized
objects, such a lump of clay, and he explicitly rejects
reducing substances to subatomic particle physics.
Already this is a red flag, since it is far from clear
that his substances are even
\emph{physically} possible;
commonsense experience with medium-sized objects is
notoriously unreliable when it comes to the peculiar
features of the quantum world such as the Heisenberg
uncertainty principle.  But let us give Nowacki the
benefit of the doubt on this point and consider the
thought experiment that he offers as an analogue of
Craig's infinite library (which we have not discussed,
but which is similar to Hilbert’s hotel):

\begin{quote}
A \emph{hyperlump} is an actually infinite lump of clay
that is composed of a denumerably infinite quantity of
different colored handfuls of clay that have been firmly
pressed together.  $\ldots\,$  The same operations Craig
performs with his actually infinite library have analogs
in the hyperlump thought experiment.  Thus, instead of a
library visitor removing books from the shelves, an artist
might approach the hyperlump and remove handfuls of clay.
The same difficulties apply as well.  For instance, it
would not be possible to add a new, numbered handful of
clay to the surface of the hyperlump: All available numbers
have already been used up in numbering the various handfuls
of clay that constitute the hyperlump.  Again, removing
handfuls of clay from the hyperlump will result in
counterintuitive absurdities.  Removing a denumerably
infinite number of handfuls of clay from the hyperlump
yields varying results: Employing one method results in
an infinite quantity of clay remaining; employing another
method wipes out the supply of clay almost entirely.

What the hyperlump example does allow us to do, however,
is bring out a family of features implicit in Craig's
thought experiment that have not been fully explained
thus far.  The nub of what is at issue is this: since
the hyperlump is a substantial body, it must have a
surface and hence possess some particular shape.  This
is because all substantial bodies have surfaces and the
surface of a body determines its shape.  However, no
particular shape we could mention is consistent with
the hyperlump.  Insofar as all substantially possible
substantial bodies necessarily have a different shape,
it follows that the hyperlump is substantially impossible.
\end{quote}

Confronted with Nowacki's claims that a hyperlump is
metaphysically impossible, we can raise exactly the same
questions that we raised before: Why are we compelled to
uphold these claims with such firmness that
we must declare any alternative
theory impossible?  To this, Nowacki does have a possible
rejoinder that was not available to Craig; namely,
if challenged to explain why a substance must have
a particular shape, he can
insist that \emph{that is just how substances are}.
But now Nowacki has set himself an even more challenging
task, which is not only to develop all of physics---quantum
field theory, general relativity, the whole shabang---on
the basis of substances, but to show, on top of all that,
that it is \emph{impossible} to develop physics
without adopting the metaphysical premises of substances.
Needless to say, Nowacki has not put in anywhere near
the amount of effort required to accomplish such a
gargantuan task.

Even if we bend over backward to cut Nowacki some slack,
and grant his theory of substances, there is still a
giant loophole in his argument.
To arrive at the conclusion that an actual infinite cannot
exist, Nowacki must give an AI-completeness proof for
hyperlumps.  That is, he has to deduce the impossibility
of an actual infinite from the impossibility of a hyperlump.
Since a hyperlump has lots of special
properties, this seems even more challenging than proving
AI-completeness for Hilbert's hotel.
How would one show that (as an example)
the possible existence of an
actual infinity of electrons implies the possible
existence of a hyperlump?  Nowacki gives us no
hint as to how such a reduction might go.

\section{An Argument from Ockham's Razor}

At this point in the discussion,
the reader might be persuaded that
we cannot rule out a priori the physical existence
of an actual infinite, but might still feel uneasy
about actively postulating an actual infinite in
the physical world.  Can we not argue as follows?
Even if it is convenient to introduce infinities
into our mathematical
models of the physical world, when it comes to making
testable predictions, we have to make
\emph{finitely} testable predictions,
for the mundane reason that we human beings are
finite creatures with finite resources at our disposal.
For example, suppose for a moment
that Hilbert's hotel exists out there
and suppose that we stumble upon it.
Even if Hilbert's hotel is actually infinite,
any prediction about Hibert's hotel that we can test
can involve only a finite portion of it.
In fact, something more is true;
if current physical theories are correct,
then there is a limit to the size
of the observable universe and so there is only a finite
portion of Hilbert's hotel that we could ever observe,
even in principle.
A finite number of observations can always be explained
by a finitary theory.
So why hypothesize an infinite wing of Hilbert's hotel,
or infinitely many invisible pages of Benardete's book,
when we can perfectly well explain all our observations with
a finitary theory?
Ockham's razor, which in one popular form states
that entities should not be multiplied without necessity,
would seem to tell us to avoid
postulating the existence of an actual infinite.
Yes, maybe postulating an actual infinite is not
as absurd and verboten as Craig would have us believe,
but if we would never postulate it anyway,
what difference does it make?

Formulating physical theories that avoid infinities
is something that some physicists do attempt from
time to time.
For example, Max Tegmark~\cite{tegmark} has publicly
argued that infinity is not needed in physics,
and that we should seek to get rid of it:

\begin{quote}
Our best computer simulations, accurately describing
everything from the formation of galaxies to tomorrow's
weather to the masses of elementary particles, use only
finite computer resources by treating everything as finite.
So if we can do without infinity to figure out what happens
next, surely nature can, too---in a way that's more deep
and elegant than the hacks we use for our computer
simulations. Our challenge as physicists is to discover
this elegant way and the infinity-free equations describing
it---the true laws of physics. 
\end{quote}

Just to be clear, I certainly do not insist that
actual infinities be postulated in physical theories.
Perhaps one day, Tegmark's dream will be realized,
and our best physical theories will avoid actual
infinities of all kinds.
However, there are two points I would like to emphasize.

First of all, it is by no means clear that Ockham's
razor requires us to excise the infinite from our theories.
Another version of Ockham's razor appeals to
\emph{simplicity} rather than multiplication of entities.
Sometimes, introducing an actual infinity into our
mathematical model yields a simpler theory than a theory
with some explicit finite bound.  Indeed, most of our best
physical theories today do employ actual infinities inside
the mathematical model.  As we alluded to earlier,
in general relativity, spacetime is modeled as a
manifold, which is traditionally thought of as an infinite
set of points.  Now, in principle, we could construct some
finitary approximation to a manifold and reconstruct all
the calculations needed for experimental predictions without
appealing to infinity.  However, in practice, physicists
and mathematicians often introduce infinities because doing
so \emph{simplifies the theory and the
calculations}.  Finite models are not always simpler,
so even if we accept Ockham's razor,
it remains an open question whether the best and simplest
physical theories will involve postulating
actual infinities.
It is therefore important to keep an open mind,
and not rule out infinities a priori,
as Craig wants us to do.

Secondly, any argument for finitary physical
theories on the grounds of our own finitude
is at best an \emph{epistemological} argument
and not a \emph{metaphysical} one.
That is, such an argument, at best, is that we have no
epistemological warrant for postulating an actual
infinite, not that we have a positive argument for the
metaphysical impossibility of an actual infinite.
But for the kalam cosmological argument, Craig needs
the actual infinite to be metaphysically impossible;
he does not just want to claim that
\emph{even if} the universe truly had no beginning then
we would not be able to know it
(or at least, we would never be able to have
confidence in a physical theory that said
that the universe had no beginning).
So even if Ockham's razor drives us to finitary theories,
it does not salvage the conclusion that Craig wants.

\section{Concluding Remarks}

As far as I have been able to tell, all the ``arguments''
proposed by supporters of the kalam cosmological argument
that an actual infinite cannot exist
in the physical world amount to nothing
more than an instinctive repugnance 
that, at bottom, is simply a metaphysical prejudice that is
best abandoned, along with all the other failed metaphysical
prejudices of the past.  But even if I am right about this,
I would like to end this essay not on a negative note, but
with two positive suggestions.

The first suggestion is that working scientists and
mathematicians take time to write down the philosophical
principles that they use to guide their own research.
These principles often go unstated, but they are valuable
and important.  Writing them down helps disseminate these
principles to more than the few privileged students
and colleagues who
have firsthand contact with successful senior researchers.
What I have been calling the prudence principle is,
I believe, widely accepted, but it is the sort of
thing that is usually not written down formally.
As I said earlier, it is not the kind of principle that
most philosophers have been trained to accept as valid,
so they are not going to appreciate the
hundreds of years of experience lying behind it unless
scientists take the time to communicate that distilled
wisdom.

My second suggestion is more speculative, and is
directed primarily to theists, especially those who
are tempted to try to salvage the
kalam cosmological argument
by rebutting the argument I have given in this essay.
All mathematicians have had the experience of trying
unsuccessfully to prove a conjecture, and then finally
coming to the realization that the reason for the
elusiveness of the proof is that the conjecture is false.
Perhaps the reason all attempts to show that an actual
infinite cannot exist have failed is that an actual
infinite \emph{can} exist (not necessarily that it
\emph{does} exist, merely that it can).
In fact, perhaps theists could try to argue for
the existence of God not from the \emph{impossibility}
of the actual infinite, but from its \emph{possibility}.

The argument that our ability to conceive of infinity is
evidence for the existence of God goes back at least
to Ren\'e Descartes's \emph{Meditations on First
Philosophy}. Descartes's argument is usually
categorized as an \emph{ontological} argument;
it does not rely on the physical existence
of an actual infinite.  Theists could consider
going further, and using the possible existence
of an actual infinite as an \emph{argument from
design} of the universe.  Though speculative, such
an argument for God would seem to be more promising
than desperate attempts to rescue the kalam
cosmological argument (at least, in the form that
I have presented it in this article).

\end{document}